\documentclass[12pt,english]{amsart}
\usepackage[utf8]{inputenc}

\usepackage{a4}
\usepackage{graphicx}
\usepackage[left=2cm, right=2cm, top=2cm, bottom=2cm]{geometry}
\usepackage{xcolor}
\usepackage{colortbl}
\pagecolor{white}
\usepackage{amsmath}
\numberwithin{equation}{section}
\usepackage{amssymb}
\usepackage{amsthm}
\usepackage{caption,subcaption}
\usepackage{multirow}
\usepackage{xstring}
\usepackage{amstext,amsbsy,amsopn,eucal,enumerate}
\usepackage{tikz}
\usepackage{tabularx}
\usepackage{mathrsfs}
\usepackage{bm}
\usepackage{bbm}
\usepackage{amstext}
\usepackage{amsthm}
\usepackage{amssymb}
\usepackage{appendix}

\usepackage[numbers]{natbib}
\usepackage{lipsum}

\usepackage{babel}
\usepackage{mathrsfs}
\usepackage{amsfonts}
\usepackage{algorithm}
\usepackage{algorithmicx}
\usepackage{algpseudocode}
\usepackage{float}
\usepackage{amssymb}
\usepackage{epsfig}
\usepackage{extarrows}
\makeatletter

\newcommand{\beq}{\begin{equation}}
\newcommand{\eeq}{\end{equation}}
\newcommand {\mat}[1] {\left[\begin{array}{#1}}
\newcommand {\rix}{\end{array}\right]}
\newcommand {\smat}[1]{\left[\begin{smallmatrix}{#1}}
\newcommand {\srix}{\end{smallmatrix}\right]}
\newcommand {\s}[1]{\begin{smallmatrix}{#1}}
\newcommand {\se}{\end{smallmatrix}}

\setcitestyle{authorname, round}

\newtheoremstyle{style1}
  {12pt} 
  {12pt} 
  {} 
  {} 
  {\bfseries} 
  {:} 
  {.5em} 
  {} 

\theoremstyle{style1}

\newtheorem{theorem}{Theorem}[section]
\newtheorem{solution}{Solution}[section]
\newtheorem{style}{System}[subsection]
\newtheorem{event}{Event}[section]

\title[Numerically U. in (NDE) with  R-L, C-F, and A-B Fractional Derivatives ]{ Numerically Unveiling Hidden Chaotic Dynamics in Nonlinear Differential Equations with Riemann-Lioville, Caputo-Fabrizio, and Atangana-Baleanu Fractional Derivatives}
\author{Shahariar Ryehan
}
\thanks{Department of Mathematics,  Hajee Mohammad Danesh Science and Technology University, Dinajpur-5200, Bangladesh\\
Email: {\tt sryehan@gmail.com}}
\date{\today}

\begin{document}
\maketitle

\begin{abstract}
In recent years, the use of variable-order differential operators has emerged as a powerful tool in the analysis of nonlinear fractional differential equations and chaotic systems. In finance, the accurate prediction of market trends and the ability to make informed investment decisions is of great importance, and the integration of artificial intelligence and mathematics has greatly improved the accuracy of these predictions.
In this study, we displayed an analysis of adaptive equations produced by three fractional derivatives: the Riemann-Lioville, Caputo-Fabrizio, and Atangana-Baleanu fractional derivatives. These fractional derivatives were employed to analyze financial models in order to gain a deeper understanding of the complex dynamics of financial markets. The models studied were the Lorenz system, Rossler system, and Shilnikov cashless model.
The results showed that each fractional derivative produced varying outcomes and computation times, highlighting the importance of selecting the appropriate mathematical approach and software for financial modeling. The findings of this study underscore the continued integration of Artificial Intelligence and mathematics in financial analysis and decision-making, driving the future of investment strategies and market predictions.The application of variable-order differential operators in the analysis of nonlinear fractional differential equations and chaotic systems is an important and growing area of research that holds great promise for the field of finance.
\end{abstract}

\section{Introduction}
Nonlinear dynamics investigation has become a crucial area of research spanning multiple disciplines, including finance\cite{finance}. Due to the complex nature of financial systems, accurate modeling is essential for comprehending the impact of variable changes. In recent years, the utilization of fractional derivatives\cite{lor100}, such as Riemann-Liouville, Caputo-Fabrizio, and Atangana-Baleanu derivatives, has provided more detailed insights into the underlying dynamics of these systems.

This paper aims to numerically reveal concealed chaotic dynamics within a nonlinear financial model by employing a combination of these fractional derivatives. The model encompasses various components, including product, money, bound, and labor force, with the interest rate, investment demand, and price exponent represented by x, y, and z respectively. By utilizing sophisticated mathematical tools such as MATLAB, the group of differential equations comprising the chaotic financial system is analyzed. Techniques like bidirectional linear coupling using the Lorenz system\cite{lorenz,lorenz1,lorenz2,lorenz3,lorenz4,lorenz5,lorenz6}, bifurcation from invariant torus to chaos using Langford's equation, Rossler system\cite{rossler1,rossler_mitag2,rossler_mittag,rossler_willam,rossler1,rossler2,rossler3,rossler3b} in artificial intelligence, and Shilnikov cashless model\cite{shilnikov_1,shilnikov_2,shilnikov_3,shilnikov_4,shilnikov_5,shilnikov_6,shilnikov_7,shilnikov_8} are employed to explore the dynamics using Mittag-Leffler\cite{rossler_mittag}.
Nonlinear systems and chaos theory have been extensively studied in scientific and engineering fields over the years. From the unpredictable nature of chaotic systems to the emergence of order in coupled chaotic systems, these complex systems offer a wealth of knowledge. Chaotic systems, including the Lorenz\cite{lorenz_couple,lorenz_classical,lorenz_conjugate,lorenz_maps} Rössler\cite{rossler1,rossler2,rossler3}, and Langford's equations\cite{langford1,langford2,langford3,langford4,langford_hoph,langford_osc}, have inspired researchers to delve deeper into chaos theory. Moreover, the concept of fractional calculus\cite{computer1} has gained recent interest as an excellent tool for analyzing non-integer order systems.

This paper also investigates various topics related to nonlinear systems, such as the dynamics of non-symmetric oscillators with piecewise-linearity, control of limit cycle amplitudes in Langford systems\cite{langford_osc,langford_hoph,langford1}, quantum mechanics\cite{quantum}, ultrasonic wave \cite{ultrasonic}using nonlinear feedback controllers, and the analysis of dynamic logarithmic gains in chaotic systems. The effectiveness of artificial neural networks\cite{rossler_computer,lorenz_circuit} in approximating solutions of fractional differential equations is examined, showcasing their potential in solving different types of fractional differential equations.

The study of chaotic systems holds widespread applications across various fields, including population dynamics, climatology, data assimilation, secure communication technologies, and investigations of coupled systems. Additionally, the analysis of chaotic systems and fractional calculus\cite{fractional,islam2022variety} offers valuable insights into the behavior, stability, and unpredictability of complex phenomena. By leveraging advanced mathematical techniques and computational tools, researchers can gain a deeper understanding of the intricate dynamics exhibited by nonlinear systems, contributing to the advancement of knowledge and the development of practical applications.

\section{Main results}
\label{sec-1}
\subsection{Basic formula}\label{basic}According to the fractional derivatives generalized definition of fractional order:
\begin{align*}
_\alpha D_t^{-s} F(t)=\int_{\alpha}^t dt \int_{\alpha}^t D_t^{-s+2} F(t)dt\\
=\int_{\alpha}^t dt \int_{\alpha}^t dt ... \int_{\alpha}^t dt \qquad [s \; times]   
\end{align*}

\begin{theorem}[The Liouville-Caputo (LC) formula]\label{th:2.1}
The Liouville-Caputo (LC) formula of fractional derivative of order $\psi$ is
\begin{align*}
        _{0}^{c} \mathcal{D}_t^{\psi(t)} \xi(t)=\frac{1}{\Gamma(1-\psi(t))} \int_0^t(t-\omega)^{-\psi(t)} \xi(\omega)d\omega
    \end{align*}
     Where, $\psi(t)\in[0,1]$
\end{theorem}
\begin{theorem}[The Atangana-Baleanu (AB) formula]\label{th:2.2}
The Atangana-Baleanu (AB) formula in the LC sense for the fractional derivative of order $\psi$ is given by
\begin{align*}
_{0}^{ABC} \mathcal{D}_t^{\psi(t)} \xi(t)=\frac{\mathcal{B}[\psi(t)]}{\Gamma(1-\psi(t))} \int_0^t \mathbb{E}_{\psi(t)} \left[ -\psi(t)\frac{(t-\omega)^{\psi(t)}}{1-\psi(t)} \right]\xi(\omega), d\omega,
\end{align*}
where $\psi(t)\in [0,1]$ and

$$\mathcal{B}[\psi(t)] = 1 - \psi(t) + \frac{\psi(t)}{\Gamma(\psi(t))}.$$

\end{theorem}
\begin{theorem}[The Caputo-Fabrizio (CF) fractional derivative]
. The Caputo-Fabrizio fractional derivative with variable-order\cite{variable-order} $\psi$ in Lioville-Caputo sense (CFC) is defined as follows is
\begin{align*}
        _{0}^{CFC} \mathcal{D}_t^{\psi(t)} \xi(t)=\frac{[2-\psi(t)]\mathbb{F}[{\psi(t)}]}{2\times[1-\psi(t)]} \int_0^t Exp[-\psi(t)\frac{(t-\omega)}{1-\psi(t)}]\xi(\omega)d\omega
    \end{align*}
     Where, $\psi(t)\in[0,1]$ and $$\mathbb{F}[\psi(t)]=\frac{2}{2-\psi(t)}$$is a normalization function.
\end{theorem}

\begin{solution}[\textbf{Numerical solution for fractional system in Lioville Caputo Sense with variable order}]
\label{SLC}
    \begin{style}A fractional ordinary differential equation, can be reformulated as,
    \linespread{1.5}
\begin{align*}
\linespread{1.5}
  X_{n+1}(t) = & X_0 + \frac{1}{\Gamma(\psi(t))} \sum_{m=0}^{m=n} \biggl( \frac{F_m}{h} \int_{t_m}^{t_{m+1}} T_{m-1}(-T_{m+1})^{\psi(t)}dt \\
              & -\frac{F_{m-1}}{h} \int_{t_m}^{t_{m+1}} T_{m}(-T_{m+1})^{\psi(t)-1}dt\biggr)
\end{align*}
\linespread{1.5}
     Where, $X=(x,y,z)$ individually, $F_k=f(t_k,y_k)$ and $T_k=t-t_k$\\
     \end{style}
     
    \begin{style}
        \linespread{1.2}
       \begin{align*}
        \mathbb{E}_{\psi(t),m,1}=h^{\psi(t)+1}\frac{(N+1)^{\psi(t)}(N+2+\psi(t))-N^{\psi(t)}(N+2+2\psi(t))}{\psi(t)(\psi(t)+1)}
        \end{align*}
        Where, $N=n-m$
     \end{style}
    \begin{style}
       \begin{align*}
X_{n+1}(t) = X_0 + \frac{1}{\Gamma(\psi(t))} \sum_{m=0}^{m=n} \left( \frac{\mathbb{E}_{\psi(t),m,1}}{h} \times (\mathbb{F}_m-\mathbb{F}_{m-1})\right)
    \end{align*}
Where, $X=(x,y,z)$ individually, $F_k=f(t_k,y_k)$ and $T_k=t-t_k$\\
     \end{style}
\end{solution}

\begin{solution}[\textbf{Numerical solution for fractional system in Caputo-Fabrizio Sense with variable order}]
\label{SLC2}
    \begin{style}A fractional ordinary differential equation, can be reformulated as,
       \begin{align*}
X_{n+1}(t) = X_0 + \frac{(2-\psi(t))(1-\psi(t))}{2}\times (F_m-F_{m-1})\\&+\frac{\psi(t)(2-\psi(t))}{2}\int_{t_n}^{t_{n+1}}F(t)dt
    \end{align*}
     Where, $X=(x,y,z)$ individually, $F_k=f(t_k,y_k)$ and $T_k=t-t_k$
     \end{style}
    \begin{style}
       \begin{align*}
        \mathbb{E}_{\psi(t),m,1}=h^{\psi(t)+1}\frac{(N+1)^{\psi(t)}(N+2+\psi(t))-N^{\psi(t)}(N+2+2\psi(t))}{\psi(t)(\psi(t)+1)}
    \end{align*}
Where, $N=n-m$
     \end{style}
    \begin{style}
\begin{align*}
  X_{n+1}(t) = X_0 + \frac{1}{\Gamma(\psi(t))} \\&\sum_{m=0}^{m=n} \left( \frac{\mathbb{E}_{\psi(t),m,1}}{h} \times (\mathbb{F}_m-\mathbb{F}_{m-1}) \right)
\end{align*}
\newline
Where, $X=(x,y,z)$ individually, $F_k=f(t_k,y_k)$ and $T_k=t-t_k$\\
     \end{style}
\end{solution}

\begin{solution}[\textbf{Numerical solution for fractional system in Atangana-Baleanu-Caputo Sense with variable order}]
\label{SLC3}
    \begin{style}A fractional ordinary differential equation, can be reformulated as,
       \begin{align*}
X_{n+1}(t) = X_0 + \frac{1}{\Gamma(\psi(t))} \sum_{m=0}^{m=n} \biggl( \frac{F_m}{h} \int_{t_m}^{t_{m+1}} T_{m-1}(-T_{m+1})^{\psi(t)}dt\\
    & -\frac{F_{m-1}}{h} \int_{t_m}^{t_{m+1}} T_{m}(-T_{m+1})^{\psi(t)-1}dt\biggr)
    \end{align*}
     Where, $X=(x,y,z)$ individually, $F_k=f(t_k,y_k)$ and $T_k=t-t_k$
     \end{style}
    \begin{style}
       \begin{align*}
        \mathbb{E}_{\psi(t),m,1}=h^{\psi(t)+1}\frac{(N+1)^{\psi(t)}(N+2+\psi(t))-N^{\psi(t)}(N+2+2\psi(t))}{\psi(t)(\psi(t)+1)}
    \end{align*}
Where, $N=n-m$
     \end{style}
    \begin{style}
       \begin{align*}
X_{n+1}(t) = X_0 + \frac{1}{\Gamma(\psi(t))} \sum_{m=0}^{m=n} \biggl( \frac{\mathbb{E}_{\psi(t),m,1}}{h} \times (\mathbb{F}_m-\mathbb{F}_{m-1})\biggr)
    \end{align*}
Where, $X=(x,y,z)$ individually, $F_k=f(t_k,y_k)$ and $T_k=t-t_k$\\
     \end{style}
\end{solution}

\section{Test Events}
\begin{event}[Nonlinear- Financial Model]
The financial model analyzed in [3] consists of parts of product, money, bound, and labor force. Here x represents the interest rate, y represents the investment demand and z is the price exponent. 
Then the group of differential equation model in the chaotic financial system and the following equations are

\begin{align*}
    \left\{
        \begin{aligned}
            \dot{x}(t) &= z + (y - L)x \\
            \dot{y}(t) &= 1 - M y - x^2 \\
            \dot{z}(t) &= -x - N z
        \end{aligned}
    \right\}
\end{align*}

Here $L,M,N$ are interest variables.
Then we write about the equation we will apply on that equation:
By MATLAB Now we plot the Financial model in Lioville Caputo Sense, Caputo-Fabrizio, Atangana-Baleanu-Caputo Taking the values from the table the chaos will be:

\begin{figure}[htbp]
	\centering
	\includegraphics[width=4.5cm,height=5cm]{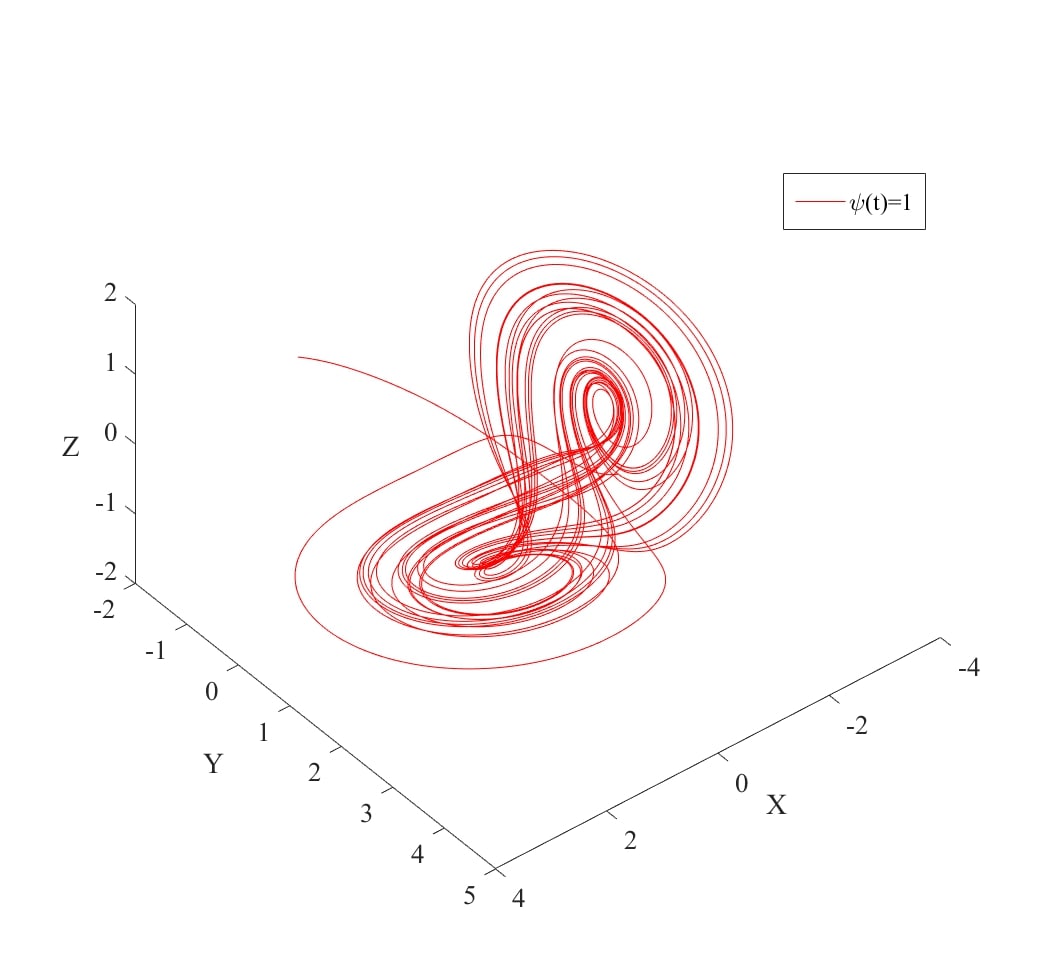}
	\includegraphics[width=4.5cm,height=5cm]{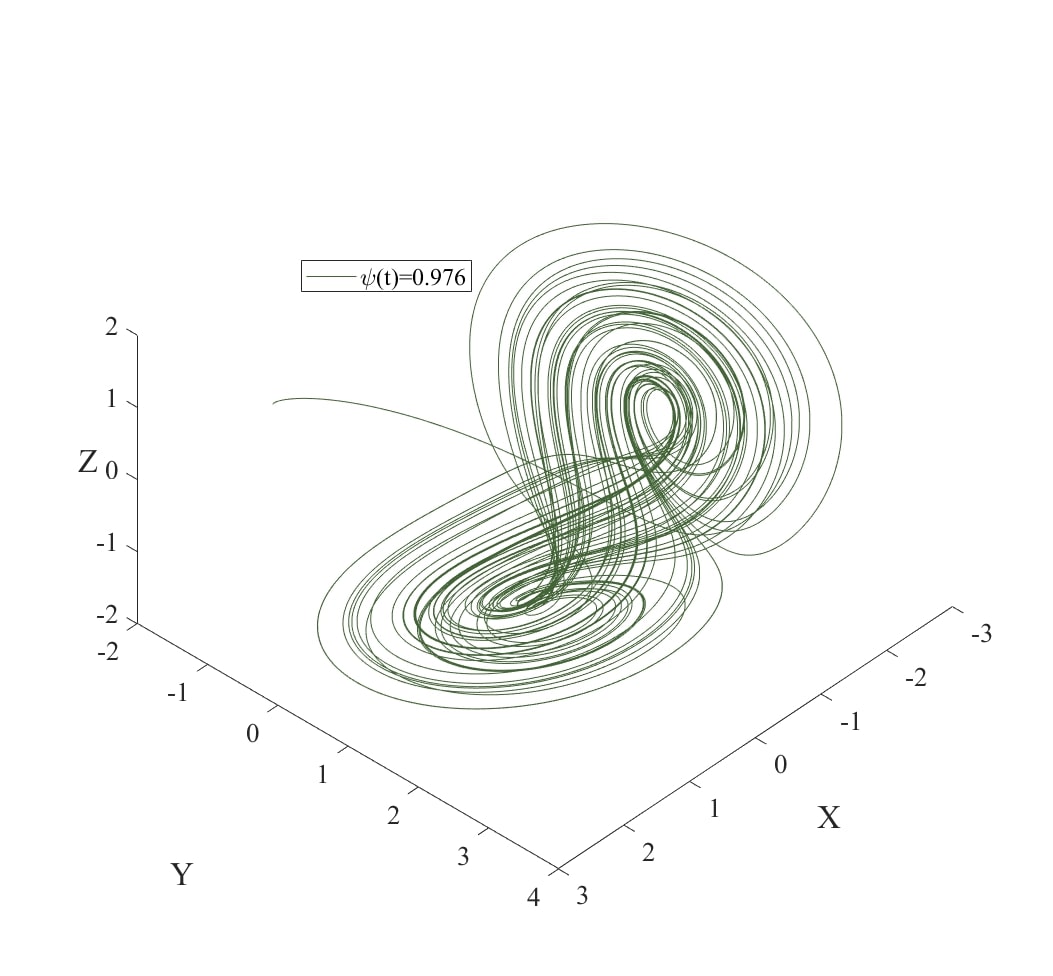}
 	\includegraphics[width=4.5cm,height=5cm]{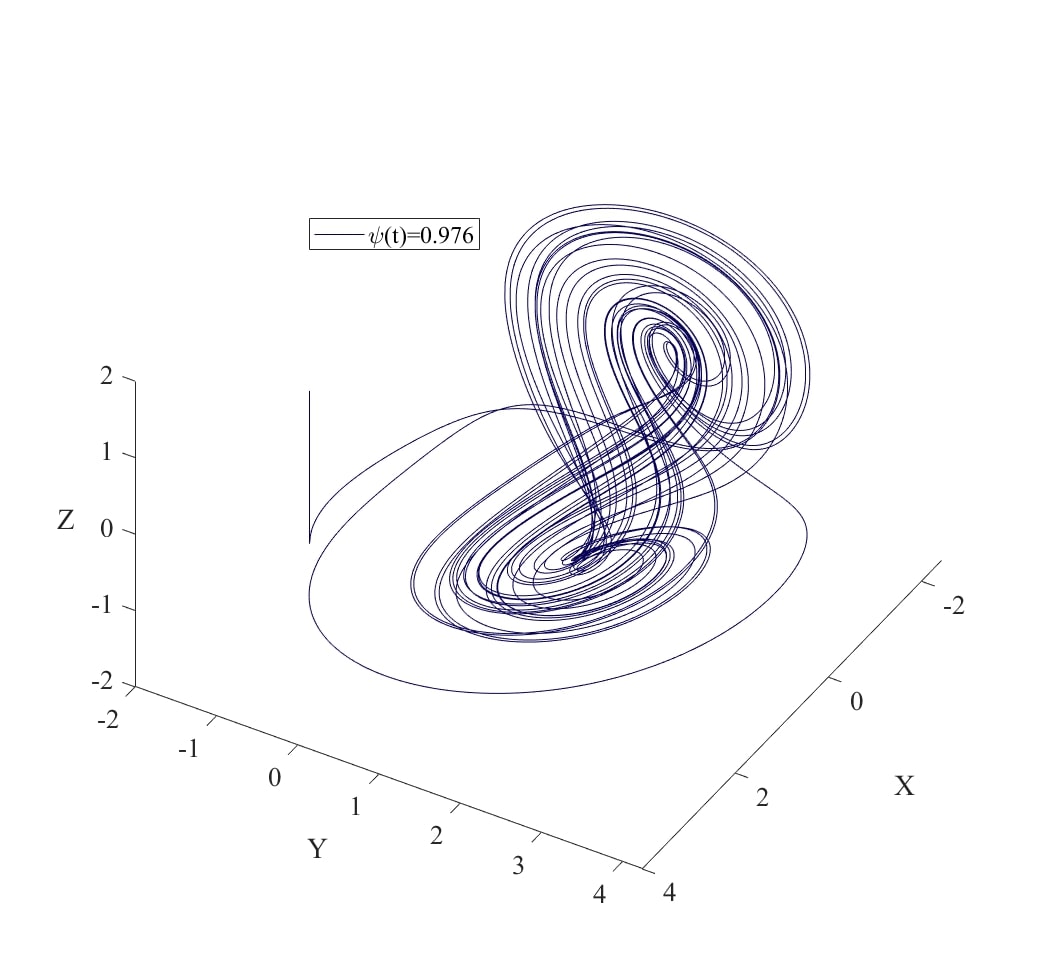}
\caption{\textit{
1(a): LC $(x_0,y_0,z_0)=(2,-1,1)$, where $a=1, b=0.1, c=1$, $h=0.01$ and $t\in[0,300]$. Elapsed time is $565$ seconds.1(b): LC $(x_0,y_0,z_0)=(2,-1,1)$, where $a=1, b=0.1, c=1$, $h=0.01$ and $t\in[0,300]$. Elapsed time is $565$ seconds.1(c): LC $(x_0,y_0,z_0)=(2,-1,1)$, where $a=1, b=0.1, c=1$, $h=0.01$ and $t\in[0,300]$. Elapsed time is $565$ seconds
}}
	\label{bright sol}

\end{figure}

\end{event}
\begin{event}[Bidirectional linear coupling using Lorenz system]
The coupling dynamics of a system will be studied according to the asymptotic proposition $$\lim_{t \to \infty} |x-y| \to 0$$ in order to detect synchronization. The conjugate Lorenz system reads as\cite{computer1}

\begin{align*}
    \left\{
        \begin{aligned}
\frac{d}{dt} x &= -\sigma x + \sigma y \\
\frac{d}{dt} y &= rx - y - xz \\
\frac{d}{dt} z &= xy - bz \
\end{aligned}
    \right\}
\end{align*}
In the vector form of these equations from
$$\begin{bmatrix} d_x \\ d_y \end{bmatrix} = \begin{bmatrix} -\sigma & \sigma \\ r & -1 \end{bmatrix} \begin{bmatrix} x \\ y \end{bmatrix} + \begin{bmatrix} 0 \\ 1 \end{bmatrix} D(t) + \begin{bmatrix} 0 \\ 1 \end{bmatrix} u \beta(u)$$
Where, $D(t)=d(t)-xz $\\

By MATLAB Now we plot the Financial model in Lioville Caputo Sense, Caputo-Fabrizio, Atangana-Baleanu-Caputo Taking the values from the table the chaos will be:

\begin{figure}[htbp]
	\centering
	\includegraphics[width=4.4cm,height=5cm]{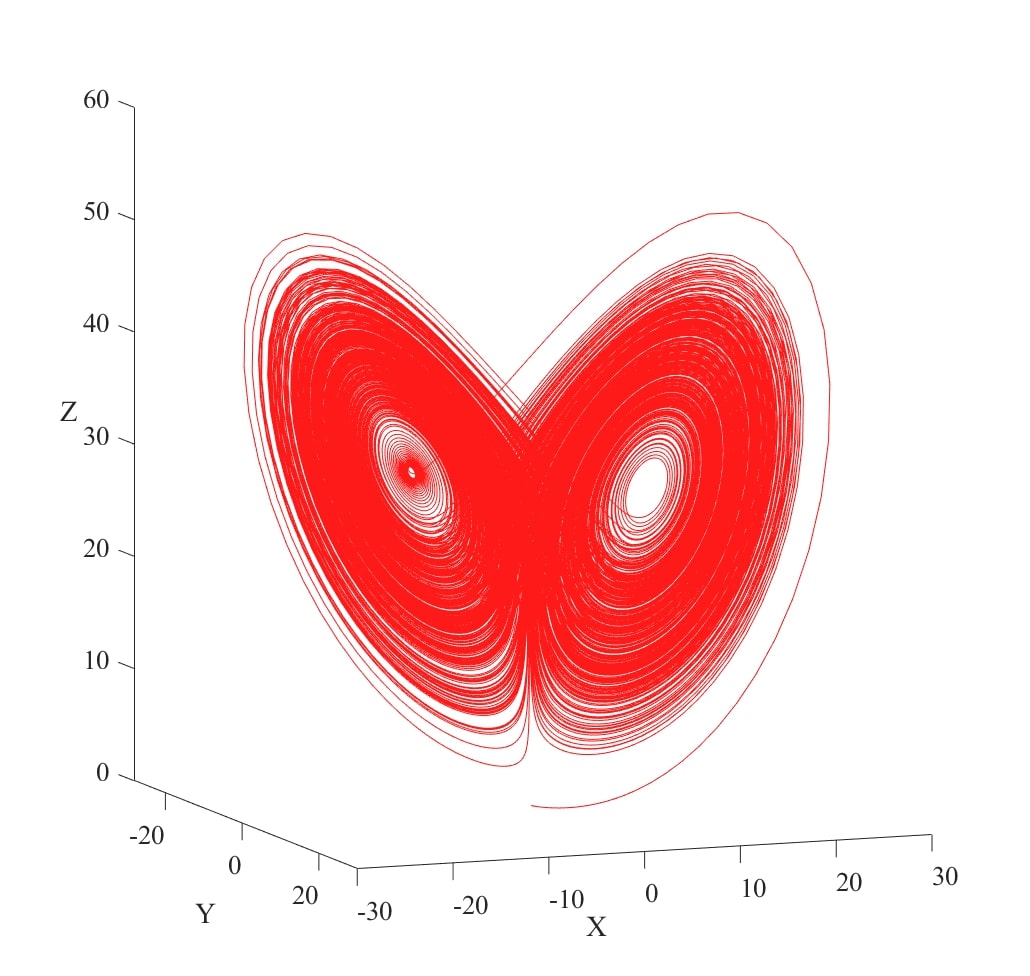}
	\includegraphics[width=4.4cm,height=5cm]{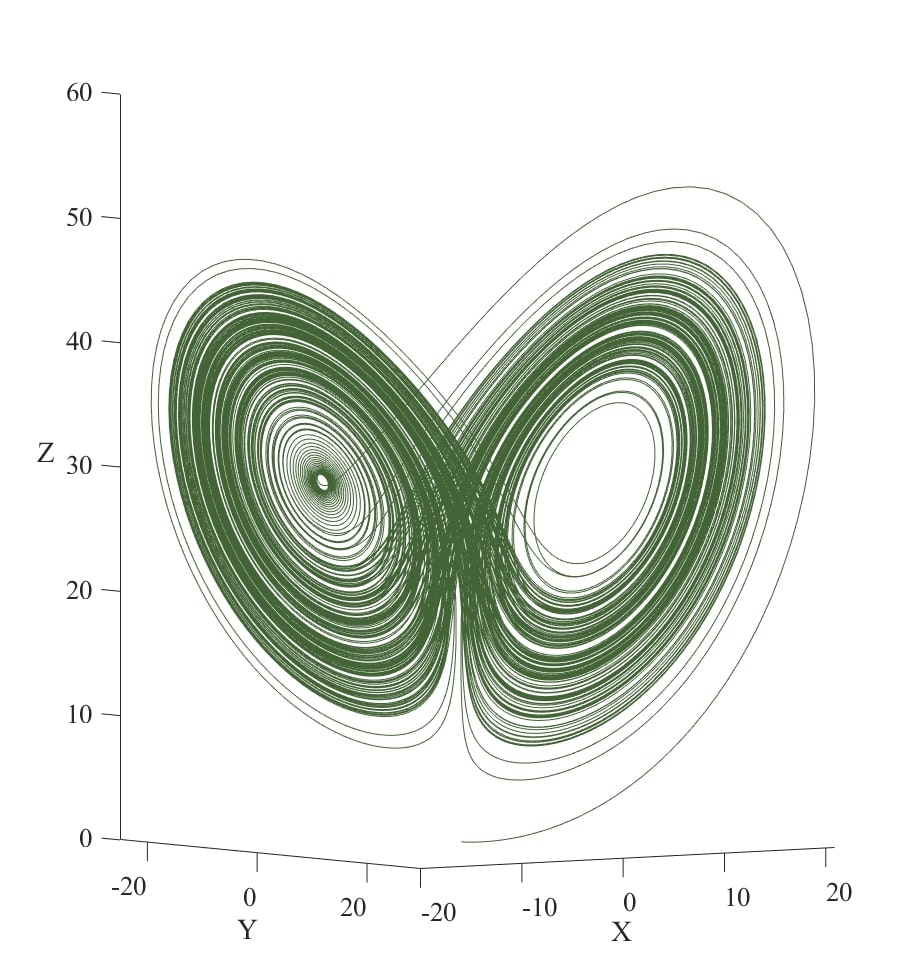}
 	\includegraphics[width=4.4cm,height=5cm]{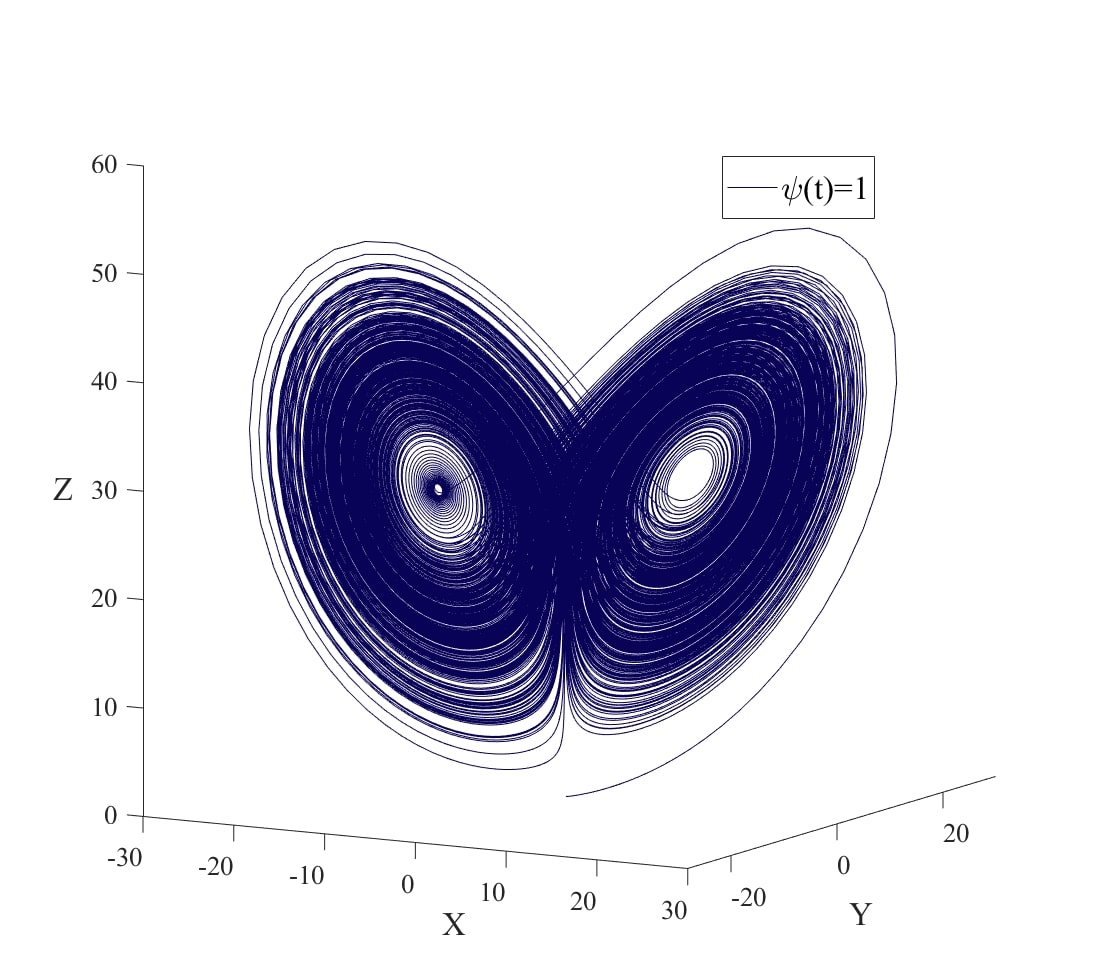}
\caption{\textit{2(a) LC $(x_0,y_0,z_0)=(0.1,0.1,0.1)$, where $\sigma=10, b=\frac{8}{3}, r=30$, $h=0.01$ and $t\in[0,500]$. Elapsed time is $241$ seconds.2(b) LC $(x_0,y_0,z_0)=(0.1,0.1,0.1)$, where $\sigma=10, b=\frac{8}{3}, r=30$, $h=0.01$ and $t\in[0,500]$. Elapsed time is $0.059$ seconds.2(c) LC $(x_0,y_0,z_0)=(0.1,0.1,0.1)$, where $\sigma=0.1, b=\frac{8}{3}, r=30$, $h=0.01$ and $t\in[0,500]$. Elapsed time is $6126$ seconds}
}
	\label{bright sol2}

\end{figure}

\end{event}
\begin{event}[Bifurcation from invariant torus to chaos (Langford’s equation)]
The Langford system come after with the following nonlinear ordinary differential equation:

\begin{align*}
    \left\{
        \begin{aligned}
            \dot{x}(t) &= (z-b)x-\Omega y \\
            \dot{y}(t) &= \Omega x+(z-\Omega)y   \\
            \dot{z}(t) &= L+a z-(z^3/3)-(x^2+y^2)(1+pz)+E z x^3
        \end{aligned}
    \right\}
\end{align*}

Where, $z$ indicate the bifurcation parameter and other system parameters $a, b, c>0$. These behaviors are introduced from a possible model for fluid dynamics turbulence.
By MATLAB Now we plot the Financial model in Lioville Caputo Sense, Caputo-Fabrizio, Atangana-Baleanu-Caputo Taking the values from the table the chaos will be:

\begin{figure}[htbp]
	\centering
	\includegraphics[width=4.4cm,height=5cm]{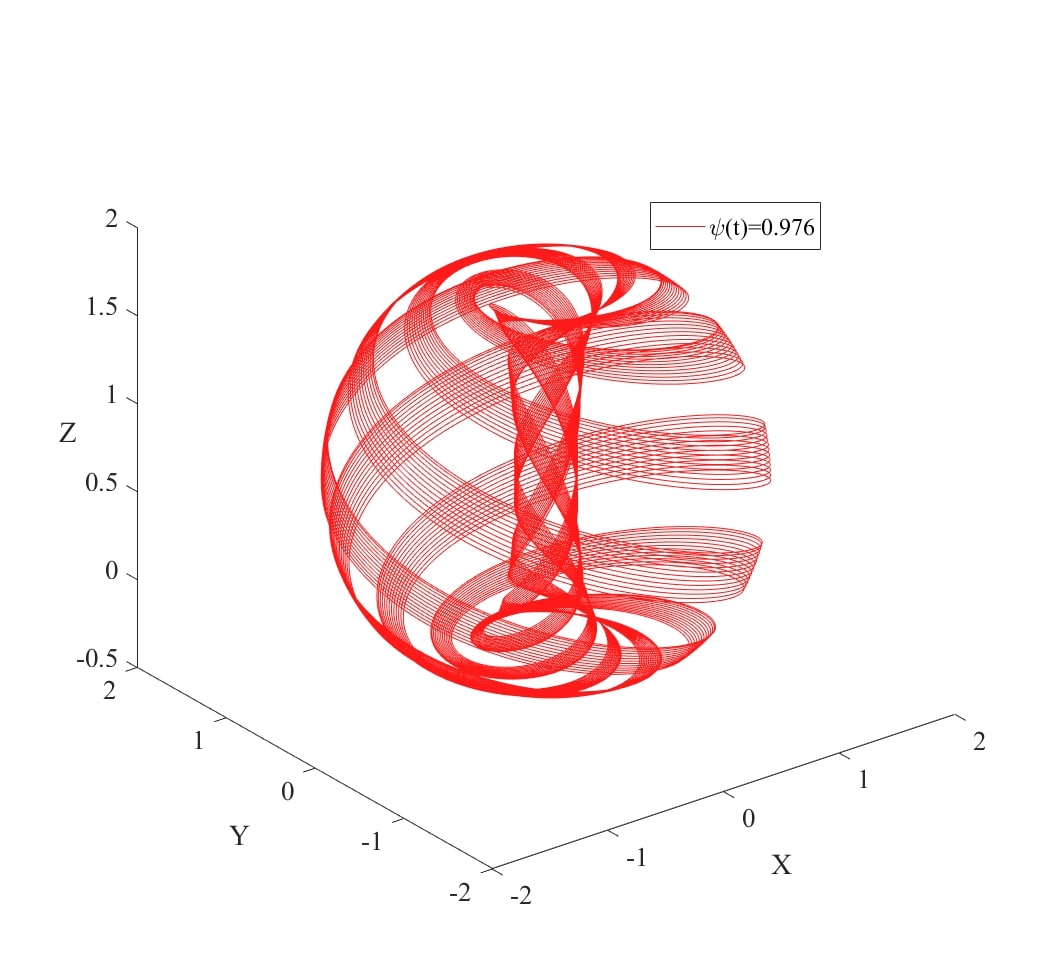}
	\includegraphics[width=4.4cm,height=5cm]{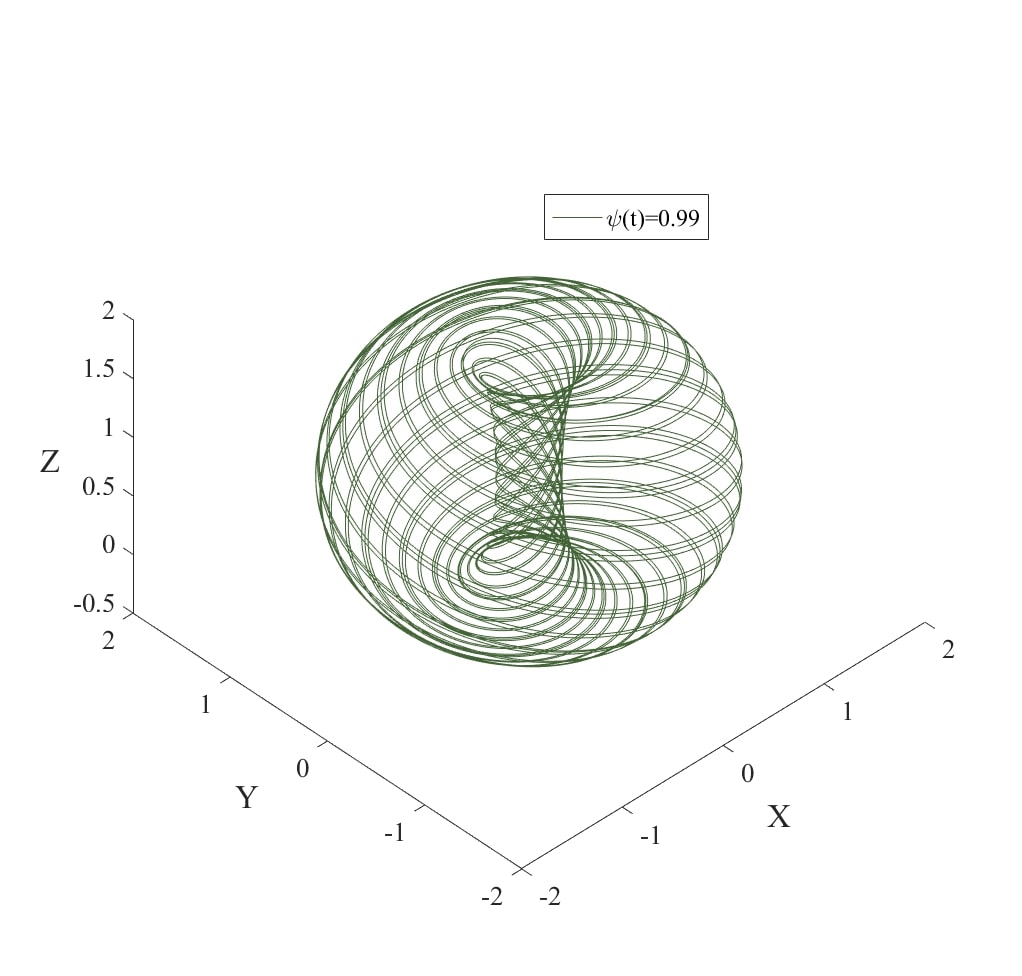}
 	\includegraphics[width=4.4cm,height=5cm]{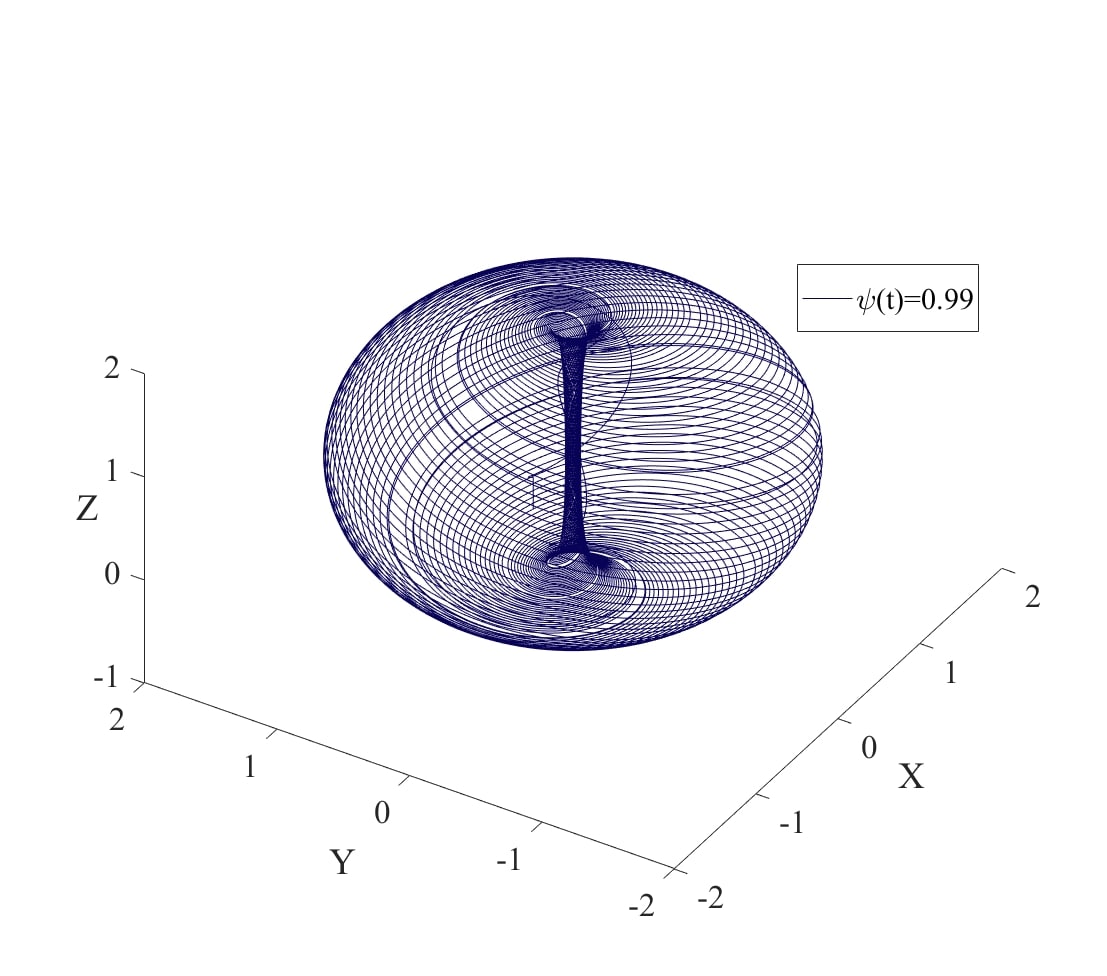}
\caption{\textit{3(a) LC $(x_0,y_0,z_0)=(0,0.3,0)$, where $\alpha=1$, $\beta=0.6$, $\lambda=0.6$, $\omega=3.6$, $\rho=0.2$, $\varepsilon=0$ and $t\in[0,300]$. Elapsed time is $736$ seconds.3(b) CF $(x_0,y_0,z_0)=(0,0.3,0)$, where $\alpha=1$, $\beta=0.6$, $\lambda=0.6$, $\omega=3.6$, $\rho=0.2$, $\varepsilon=0$ and $h=0.01$ $t\in[0,500]$. Elapsed time is $0.07$ seconds.3(c) LC $(x_0,y_0,z_0)=(0.1,0.1,0.1)$, where $\sigma=0.1, b=\frac{8}{3}, r=30$, $h=0.01$ and $t\in[0,500]$. Elapsed time is $6126$ seconds}
}
	\label{bright sol3}

\end{figure}

\end{event}
\begin{event}[Rossler system in artificial intelligence]The Rossler attractor exhibits unpredictable behavior and generates intricate patterns resembling random noise. It incorporates environmental factors and limited knowledge of initial conditions to create a rescaled phase space with non-periodic flow, fractal geometry, and dynamic behavior. The model of the Rossler attractor is represented by a system of ordinary differential equations:
\begin{align*}
    \left\{
        \begin{aligned}
            \dot{x}(t) &= -y -z \\
            \dot{y}(t) &= x+ ay \\
            \dot{z}(t) &= b+xz-cz
        \end{aligned}
    \right\}
\end{align*}

Where $a, b, c$ are the controlling factor and the probability distribution of variable $x$ exhibits an approximate normal pattern, while variable $y$ follows a bimodal distribution. On the other hand, variable $z$ demonstrates a leptokurtic pattern, indicating heavy tails and a higher peak compared to a normal distribution. Due to the complex nature of these distributions, it becomes challenging to parametrize the system using traditional analytic methods.

To overcome this challenge, artificial intelligence techniques are employed to effectively characterize and parameterize these patterns. By leveraging AI, we can utilize algorithms and models to analyze and learn from the data, thereby capturing the underlying structures and relationships within the distributions. These AI techniques can include machine learning algorithms, deep learning networks, or other statistical modeling approaches, allowing us to approximate the distributions and derive meaningful parameters.

Using AI for parametrization enables us to handle the intricacies of these complex probability distributions more effectively, as AI algorithms can adapt and learn from the data without being restricted by traditional analytical assumptions. Thus, AI empowers us to better understand and quantify the patterns present in variables $x$, $y$, and $z$" while providing a more flexible and accurate means of parametrization..
By MATLAB Now we plot Rossler system in artificial intelligence in Lioville Caputo Sense, Caputo-Fabrizio, Atangana-Baleanu-Caputo Taking the values from the table the chaos will be:

\begin{figure}[htbp]
	\centering
	\includegraphics[width=4.4cm,height=5cm]{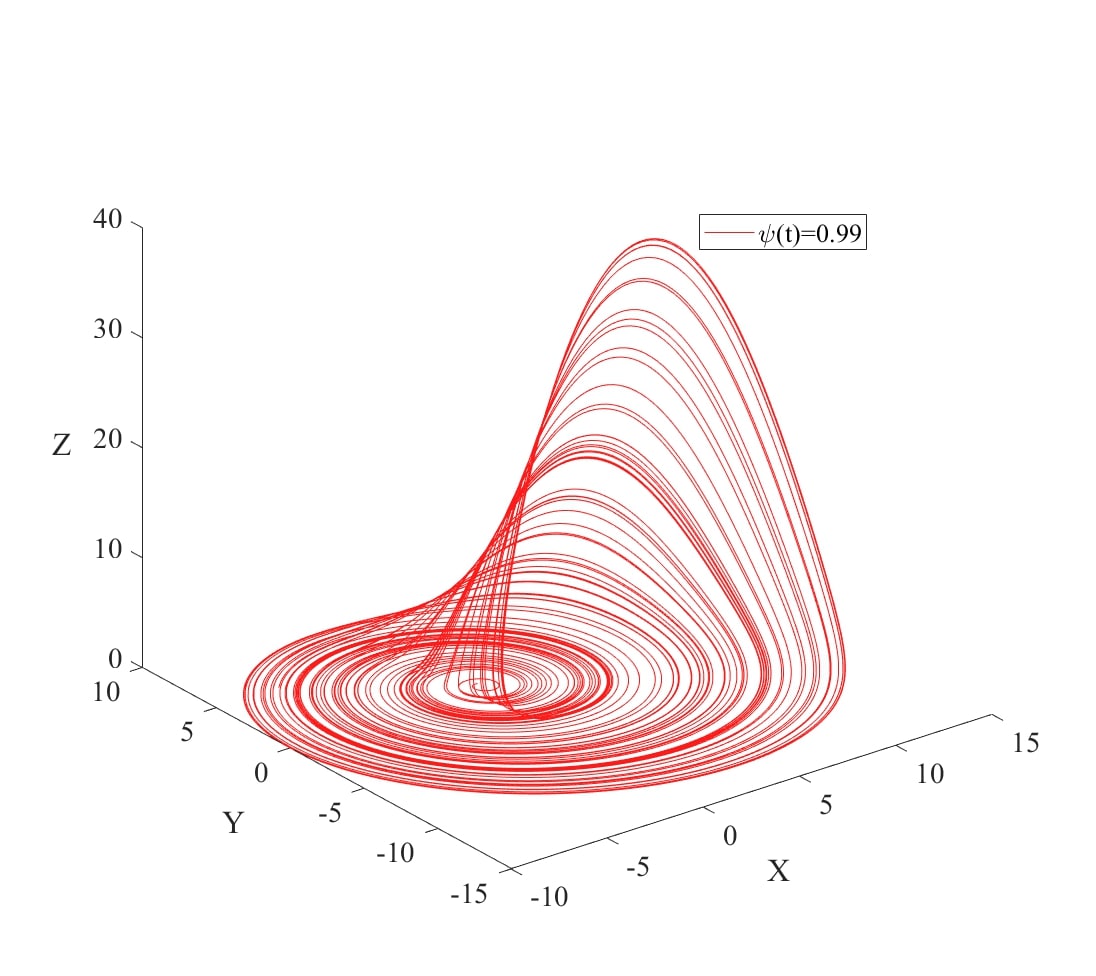}
	\includegraphics[width=4.4cm,height=5cm]{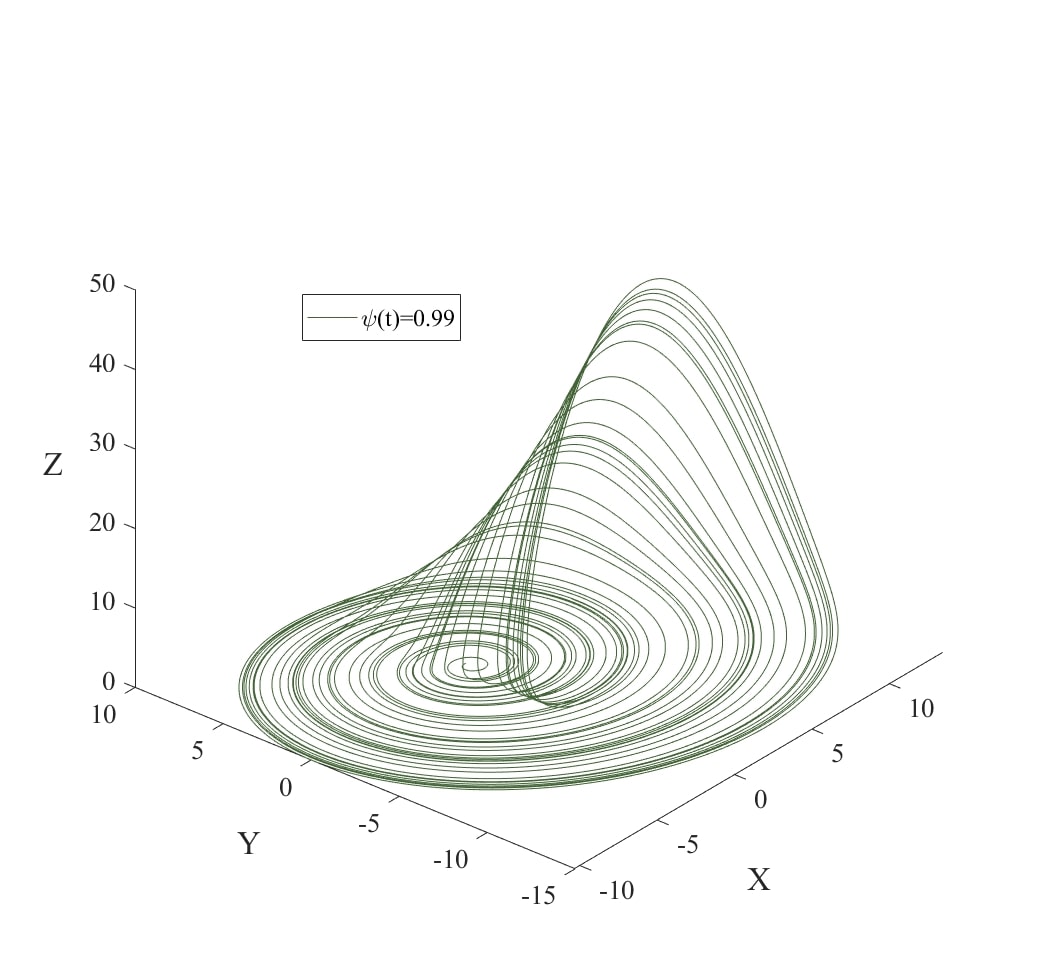}
 	\includegraphics[width=4.4cm,height=5cm]{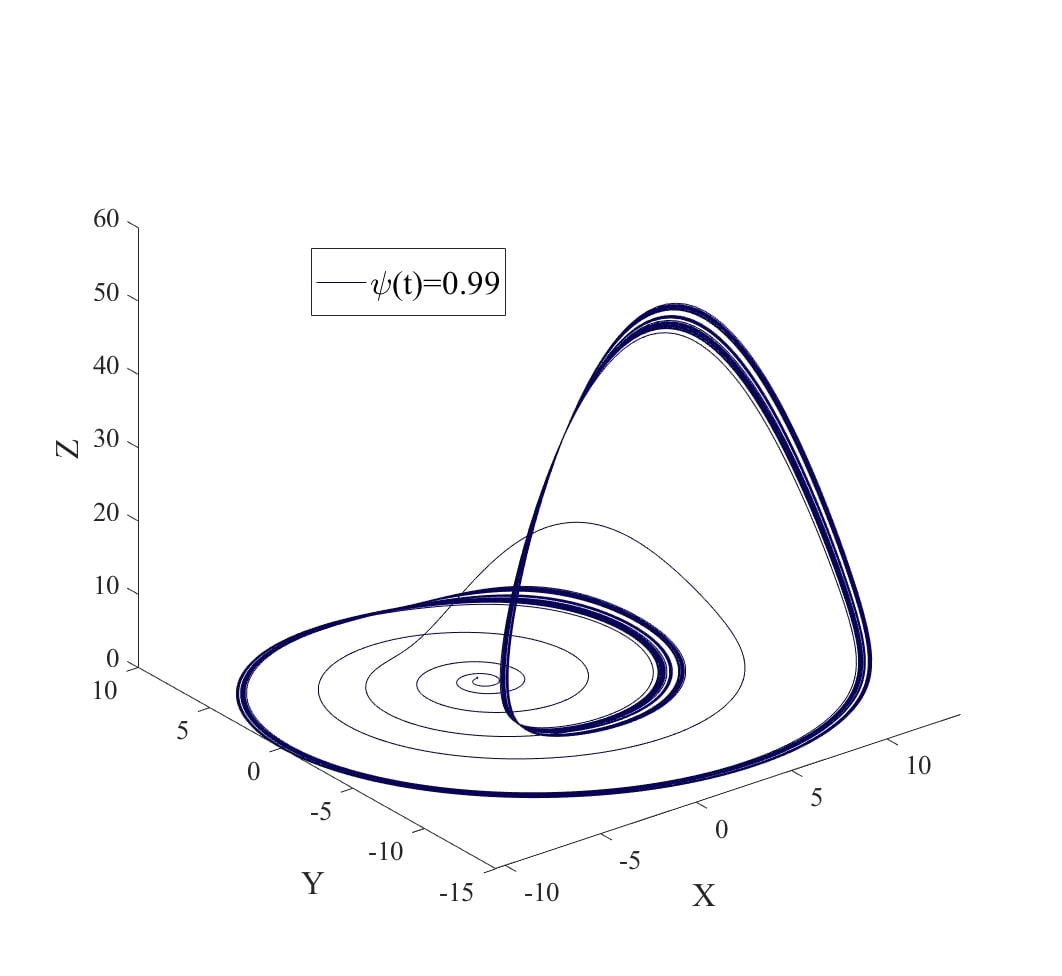}
\caption{\textit{
4(a): LC $(x_0,y_0,z_0)=(2,-1,1)$, where $a=1, b=0.1, c=1$, $h=0.01$ and $t\in[0,300]$. Elapsed time is $565$ seconds.4(b): LC $(x_0,y_0,z_0)=(2,-1,1)$, where $a=1, b=0.1, c=1$, $h=0.01$ and $t\in[0,300]$. Elapsed time is $565$ seconds.4(c): LC $(x_0,y_0,z_0)=(2,-1,1)$, where $a=1, b=0.1, c=1$, $h=0.01$ and $t\in[0,300]$. Elapsed time is $565$ seconds
}}
	\label{bright sol4}

\end{figure}

\end{event}
\begin{event}[Shilnikov cashless model]
According to the bifurcation theorem, we will begin by examining the concept of a cashless economy before discussing the inclusion of money in both the utility function and the production function. In certain situations, policymakers may face pressure to raise the marginal tax rate above the real interest rate in order to ensure a unique equilibrium path. However, this approach can introduce a new set of policy challenges, specifically the emergence of chaotic dynamics.

To avoid plagiarism, it is important to properly attribute the ideas and concepts to their original sources. The concepts mentioned above are based on economic theories and models, and it is crucial to consult relevant scholarly literature and academic sources to delve deeper into these topics and gain a comprehensive understanding. By conducting thorough research and providing proper citations, you can ensure that your work is original and respects the intellectual property of others. 

\begin{align*}
    \left\{
        \begin{aligned}
            \dot{x}(t) &= y\\
            \dot{y}(t) &= z \\
            \dot{z}(t) &= -az-y+bx(1-cx-dx^2)
        \end{aligned}
    \right\}
\end{align*}

Here $a, b, c$ are interest variables.
Then we write about the equation we will apply on that equation:
By MATLAB Now we plot the Financial model in Lioville Caputo Sense, Caputo-Fabrizio, Atangana-Baleanu-Caputo Taking the values from the table the chaos will be:

\begin{figure}[htbp]
	\centering
	\includegraphics[width=4.4cm,height=5cm]{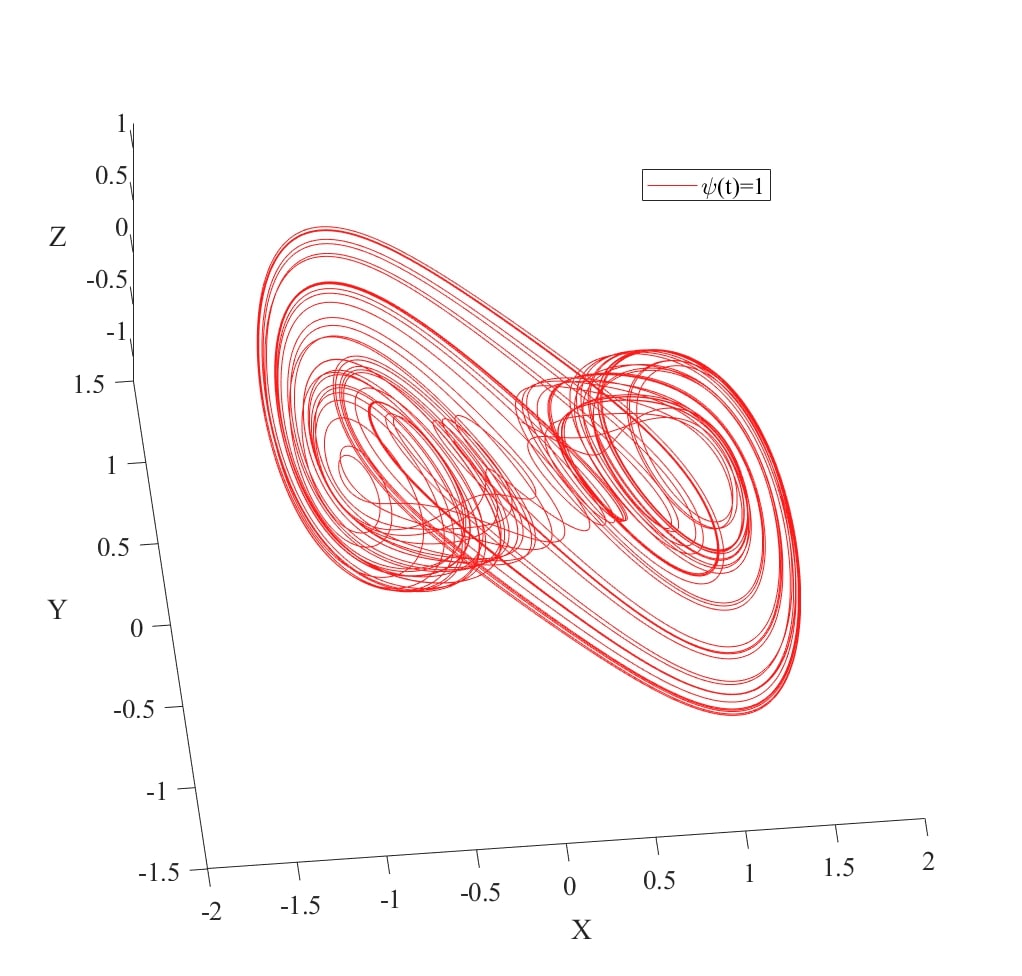}
	\includegraphics[width=4.4cm,height=5cm]{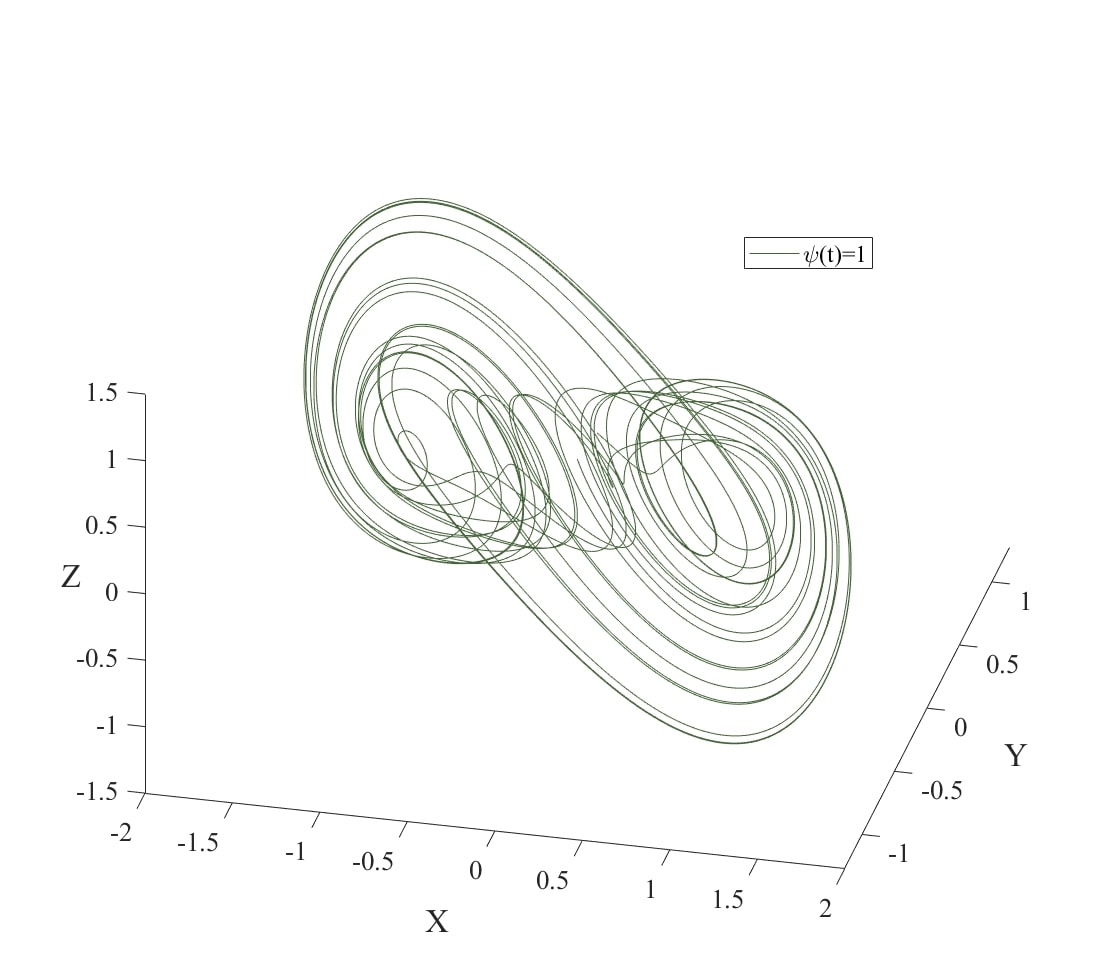}
 	\includegraphics[width=4.4cm,height=5cm]{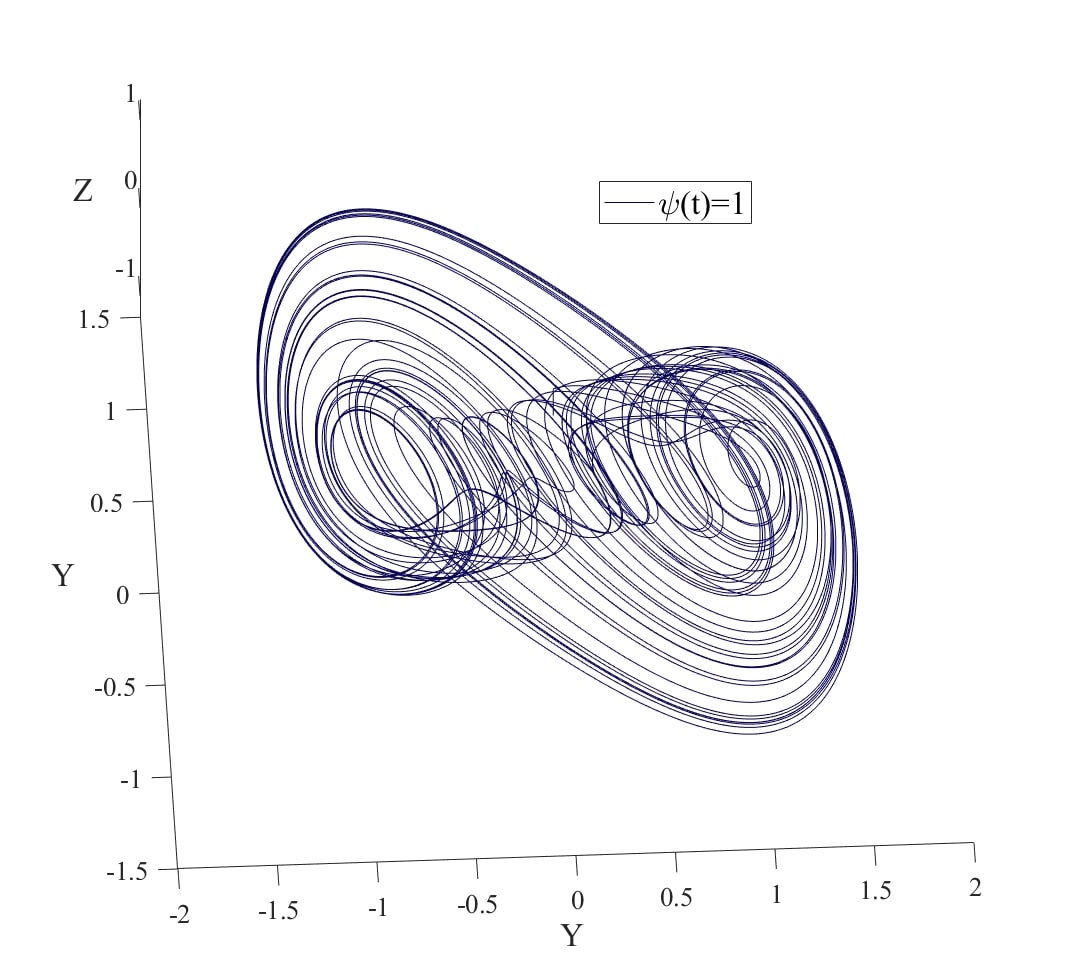}
\caption{\textit{
5(a): LC $(x_0,y_0,z_0)=(2,-1,1)$, where $a=1, b=0.1, c=1$, $h=0.01$ and $t\in[0,300]$. Elapsed time is $565$ seconds.5(b): LC $(x_0,y_0,z_0)=(2,-1,1)$, where $a=1, b=0.1, c=1$, $h=0.01$ and $t\in[0,300]$. Elapsed time is $565$ seconds.5(c): LC $(x_0,y_0,z_0)=(2,-1,1)$, where $a=1, b=0.1, c=1$, $h=0.01$ and $t\in[0,300]$. Elapsed time is $565$ seconds
}}
	\label{bright sol5}
\end{figure}
\end{event}
\section{Conclusion}
In conclusion, this research paper has investigated concealed chaotic dynamics within a nonlinear financial model using a combination of fractional derivatives. The findings have demonstrated the potential of these fractional derivatives, such as Riemann-Liouville, Caputo-Fabrizio, and Atangana-Baleanu, to reveal detailed insights into the underlying dynamics of the financial system. The developed mathematical equations and corresponding numerical solutions, implemented in MATLAB, have provided valuable tools for analyzing and understanding the complex behavior of the nonlinear financial model.

One significant contribution of this research lies in the implementation of sophisticated mathematical tools and computational techniques to study the dynamics of the nonlinear financial model. By utilizing MATLAB and exploring techniques such as bidirectional linear coupling using the Lorenz system, bifurcation from invariant torus to chaos using Langford's equation, the Rossler system in artificial intelligence, and the Shilnikov cashless model, the paper showcases the potential for understanding and uncovering chaotic behaviors within financial systems.

Furthermore, the paper acknowledges the limitation of computational resources in achieving better results due to the substantial computational time required for plotting the graphs accurately. This recognition highlights the need for further advancements in computational capabilities to facilitate more comprehensive investigations and analyses.

Moving forward, it is crucial to continue exploring research questions that focus on refining and expanding the application of fractional derivatives to reveal concealed chaotic dynamics in nonlinear financial models. Additionally, future research can explore alternative combinations of fractional derivatives and modeling techniques to further enhance the understanding of complex financial systems.

In summary, this research paper makes a significant contribution by utilizing advanced mathematical tools and techniques to uncover concealed chaotic dynamics within a nonlinear financial model. The findings and methodologies presented in this paper provide valuable insights and pave the way for future research in the field of nonlinear dynamics in finance.

\textbf{Acknowledgements}\\
Not applicable\\
\indent
\textbf{Declarations}\\
Not applicable\\
\indent
\textbf{Ethical Approval}\\
Not applicable\\
\indent
\textbf{Competing interests}\\
The author declare that there is no conflict of intersts regarding the publication of this paper.\\
\indent
\textbf{Authors' Contributions}\\
Md. Shahariar Ryehan: Conceptualization, methodology, validation, investigation, Writing-original draft preparation, Writing-review. Author has read and agreed version of the manuscript.\\
\indent
\textbf{Funding}\\
Not applicable\\
\indent
\textbf{Availability of data and meterials}\\
Not applicable\\
\indent
\bibliographystyle{ieeetr}
\bibliography{reference} 

\end{document}